\documentclass[12pt]{amsart}
\usepackage{amsmath, latexsym, amssymb, amsfonts,amscd,enumerate}
\usepackage{epstopdf}
\usepackage{graphicx}
\usepackage{fullpage}
\usepackage{cite}
\usepackage{epstopdf}
\usepackage{subfigure}

\title{Colorful Proofs of the Generating Formulas for Signed and Unsigned Stirling Numbers of the First Kind}
\author{Paul Levande}

\begin{document}
\maketitle
\begin{abstract}
We describe proofs of the standard generating formulas for unsigned and signed Stirling numbers of the first kind that follow from a natural combinatorial interpretation based on cycle-colored permutations. 
\end{abstract}
\section{Introduction}
For $\pi \in S_{n}$, let $k(\pi)$ be the number of cycles of $\pi$ when $\pi$ is uniquely written as a product of disjoint cycles (including fixed points).  Then $x^{k(\pi)}$ is the number of ways to color the cycles of $\pi$ using $x$ possible colors, and the following equations hold 
\begin{eqnarray}
\sum_{\pi \in S_{n}} x^{k(\pi)} &=& (x+n-1)(x+n-2) \cdots x \\
\sum_{\pi \in S_{n}} (-1)^{n-k(\pi)}x^{k(\pi)} &=& (x-n+1)(x-n+2) \cdots x
\end{eqnarray}
These equations are usually written as 
\begin{eqnarray}
\sum_{k} c(n, k)x^{k} &=& (x+n-1)(x+n-2) \cdots x \\
\sum_{k} (-1)^{n-k}c(n,k)x^{k} &=& (x-n+1)(x-n+2) \cdots x
\end{eqnarray}
where $c(n, k)$ is the number of permutations of $[n]$ with $k$ cycles.  The numbers $c(n, k)$ as $n, k$ vary are the \emph{unsigned Stirling numbers of the first kind}, with the numbers $(-1)^{k}c(n, k)$ the \emph{signed Stirling numbers of the first kind}, and equations $(1)$ and $(2)$ can therefore be interpreted as the generating functions of the signed and unsigned Stirling numbers of the first kind.  Stanley~\cite{stanbook} gives three different proofs of the generating formula for unsigned Stirling numbers of the first kind and proves the generating formula for signed Stirling numbers of the first kind as an algebraic consequence.  We will provide separate bijective proofs of both generating formulas, based on the first interpretations given in this paper, where both formulas count the ordered pairs $(\pi, \mu)$, where $\pi \in S_{n}$ and $\mu$ is a coloring of the cycles of $\pi$ using $x$ colors, or the number of \emph{cycle-colored permutations of $[n]$}.   
\section{The Generating Formula for Unsigned Stirling Numbers}
To prove \[
\sum_{\pi \in S_{n}} x^{k(\pi)} = (x+n-1)(x+n-2) \cdots x
\]
it suffices to define a bijection between cycle-colored permutations of $[n]$, i.e., between ordered pairs $(\pi, \mu)$, where $\pi \in S_{n}$ and $\mu$ is a coloring of the cycles of $\pi$ using $x$ colors, and sequences $a_{1} a_{2} \ldots a_{n}$ with no repeated elements such that $\left\{a_{1}, a_{2}, \ldots, a_{n} \right\} \subset X$, where $X$ is the set consisting of $2, 3, \ldots, n$ and $x$ letters $r, b, g, \ldots$.    (The proof below does not require the letters be ordered).  To define such a bijection, we first introduce a new set, the set of \emph{cycle-colored permutation graphs}.  
\subsection{Cycle-colored Permutation Graphs}
Given $A \subset X -\left\{2, 3, \ldots, n \right\}$ a subset of the set of $x$ letters $r, b, g, \ldots$, let $Y$ be a \emph{directed} graph on $[n] \cup A$.  $Y$ is a \emph{cycle-colored permutation graph} if and only if:
\begin{itemize}
	\item For all $i \in [n]$, there is a unique $j \in [n] \cup A$ such that $(i, j)$ is an edge of $Y$
	\item For all $1 \neq i \in [n] \cup A$, there is a unique $k \in [n]$ such that $(k, i)$ is an edge of $Y$
	\item If $i \in [n] \cup A$, $(i, 1)$ is not an edge of $Y$
	\item if $i \in A$, $(s, i)$ is not an edge of $Y$ for any $s \in A$
\end{itemize}
Note that there is a bijection between cycle-colored permutation graphs on $[n] \cup A$ and sequences $a_{1} a_{2} \ldots a_{n}$ such that $A \subset \left\{a_{1}, a_{2}, \ldots, a_{n} \right\} \subset \left\{2, 3, \ldots, n \right\} \cup A$ (i.e., such that all the letters of $A$ are used in $a_{1} a_{2} \ldots a_{n}$): Given a cycle-colored permutation graph $Y$ on $[n] \cup A$, define a sequence $a_{1} a_{2} \ldots a_{n}$, where $a_{i} \in [n] \cup A$ is the unique element such that $(i, a_{i})$ is an edge of $Y$.  Conversely, given a sequence $a_{1} a_{2} \ldots a_{n}$ such that $\left\{a_{1}, a_{2}, \ldots, a_{n} \right\} \subset \left\{2, 3, \ldots, n \right\} \cup A$, define a cycle-colored permutation graph $Y$ on $[n] \cup A$ by including precisely the edges $(i, a_{i})$ for $1 \leq i \leq n$.  Using this simple bijection, to prove the generating function for unsigned Stirling numbers, it suffices to define a bijection between cycle-colored permutations $(\pi, \mu)$, where $\pi \in S_{n}$ and $\mu$ colors the cycles of $\pi$ using the colors $A$, and cycle-colored permutation graphs $Y$ on $[n] \cup A$.  
\subsection{From Cycle-Colored Permutations to Cycle-Colored Permutation Graphs}
Given a cycle-colored permutation $(\pi, \mu)$, where $\pi \in S_{n}$ and $\mu$ colors the cycles of $\pi$ using the colors $A$, define a directed graph $Y$ on $[n] \cup A$ as follows:
\begin{itemize}
	\item For all $s \in A$, let $B^{s} \subset [n]$ be the set of elements of $[n]$ in cycles of $\pi$ colored using $s$ under $\mu$, with $c_{1}^{s} < c_{2}^{s} < \ldots < c_{M}^{s}$ the unique increasing rearrangement of $B^{s}$, and $M = |B^{s}|$. 
	\item For all $s$ such that $1 \notin B^{s}$, define the sequence $b_{1}^{s}b_{2}^{s} \ldots b_{M}^{s}$ by $\pi(c_{i}^{s}) = b_{i}^{s}$.  Alternately, $b_{1}^{s}b_{2}^{s} \ldots b_{M}^{s}$ is the permutation $\pi|_{B^{s}}: B^{s} \rightarrow B^{s}$ written in list form.  
	\item Define a directed graph $Y$ on $[n] \cup A$ by including, for $s$ such that $1 \notin B^{s}$, precisely the directed edges $(b_{i}^{s}, b_{i+1}^{s})$ for $1 \leq i < |B^{s}|$ and $(b_{|B^{s}|}^{s}, s)$; as well as, for the unique $s$ such that $1 \in B^{s}$, the directed edges $(i, \pi(i))$ for all $i \neq \pi^{-1}(1) \in B^{s}$ as well as the directed edge $(\pi^{-1}(1), s)$.    
\end{itemize}
Then $Y$ is a cycle-colored permutation graph on $[n] \cup A$.  

For all $1 \leq i \leq n$, let $a_{i}$ be the unique element of $[n] \cup A - \left\{1 \right\}$ such that $(i, a_{i})$ is an edge of $Y$.  Then $a_{1} a_{2} \ldots a_{n}$ is a sequence without repeated elements such that $A \subset \left\{a_{1}, a_{2}, \ldots, a_{n} \right\} \subset \left\{2, 3, \ldots, n \right\} \cup A$.

\subsection{Example} Let $n=11$, $A = \left\{r, b, g \right\}$, and $(\pi, \mu) = ( (1, 11, 3) (2, 9) (4) (5,7) (6, 8, 10), \mu)$, where $\mu$ is the coloring of the cycles of $\pi$ which has $(1, 11, 3)$ and $(6, 8, 10)$ colored using the color $r$; $(2, 9)$  and  $(4)$ colored using the color $b$, and $(5,7)$ colored using the color $g$.  Then $B^{r} = \left\{1, 3, 6, 8, 10, 11 \right\}$, $B^{b} = \left\{2, 4, 9 \right\}$, and $B^{g} = \left\{5, 7 \right\}$.  $1 \notin B^{b}, B^{r}$, and so $c_{1}^{b} = 2$, $c_{2}^{b} = 4$, $c_{3}^{b} = 9$; with $b_{1}^{b} = \pi(2) = 9$; $b_{2}^{b} = \pi(4) = 4$, and $b_{3}^{b} = \pi(9) = 2$; $c_{1}^{g} = 5$, $c_{2}^{g} = 7$, and so $b_{1}^{g} = 7$, $b_{2}^{g} = 5$.  Define $Y$ on $[11] \cup \left\{r, b, g \right\}$ by including precisely the edges $(9, 4)$, $(4, 2)$, $(2, b)$; $(7, 5)$, $(5, g)$; and $(1, 11)$, $(11, 3)$, $(6, 8)$, $(8, 10)$, $(10, 6)$, $(3, r)$, and $Y$ is the directed graph shown in Figure $1$, with $a_{1} a_{2} \ldots a_{11}$ = $11, b, r, 2, g, 8, 5, 10, 4, 6, 3$.  

\begin{center}
\includegraphics[scale=.5]{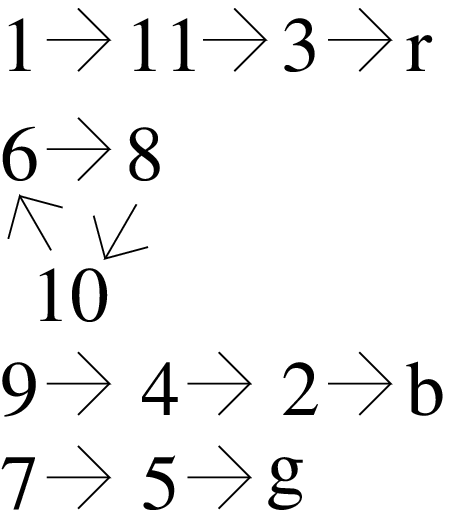}
\\Figure 1
\end{center}
\subsection{From Cycle-Colored Permutation Graphs to Cycle-colored Permutations}
Given a cycle-colored permutation graph $Y$ on $[n] \cup A$, or equivalently, given a sequence $a_{1} a_{2} \ldots a_{n}$ without repeated elements such that $ A \subset \left\{a_{1}, a_{2}, \ldots, a_{n} \right\} \subset \left\{2, 3, \ldots, n \right\} \cup A$, define a cycle-colored permutation $(\pi, \mu)$, where $\pi \in S_{n}$ and $\mu$ is a coloring of the cycles of $\pi$ using precisely the colors in $A$, as follows:
\begin{itemize}
	\item Define a permutation $\pi \in S_{n}$ as follows:  For all $s \in A$, let $B^{s}$ be the set of $i \in [n]$ such that $i$ and $s$ are connected in $Y$, with $c_{1}^{s} < c_{2}^{s} < \ldots < c_{M}^{s}$ the unique increasing rearrangement of $B^{s}$, and $M = |B^{s}|$.  Also include, in the unique $B^{s}$ such that $1$ is connected to $s$, every $i \in [n]$ unconnected to \emph{any} $s \in A$.  For $s$ such that $1 \notin B^{s}$, define a sequence $b_{1}^{s} b_{2}^{s} \ldots b_{M}^{s}$, where $b_{1}^{s} \in B^{s}$ is the element with the longest (unique) path to $s$ and $b_{1}^{s}, b_{2}^{s}, \ldots, b_{M}^{s}, s$ is that path.  Then $\pi(c_{i}^{s}) = b_{i}^{s}$.  For the unique $s$ such that $1 \in B^{s}$, for all $i \in B^{s}$, either $(i, j)$ is an edge of $Y$ where $j$ is uniquely determined, in which case $\pi(i) = j$; or $(i, s)$ is an edge of $Y$, in which case $\pi(i) = 1$,  
	\item Define a coloring $\mu$ of the cycles of $\pi$ as follows: if $i \in B^{s}$, the cycle containing $i$ is colored using $s$.
	\end{itemize}
and $(\pi, \mu)$ is a cycle-colored permutation of $[n]$.  The parallel notations used to describe the maps between cycle-colored permutations and cycle-colored permutation graphs should suffice to show that the maps invert each other, which suffices to prove the above generating function of unsigned Stirling numbers.    
\subsection{Example} Let $n=11$, $A = \left\{r, b, g\right\}$, and $a_{1} a_{2} \ldots a_{11} = 5, 6, 7, 4, 10, b, r, 3, 11, g, 9$.  Equivalently, let $Y$ be the cycle-colored permutation graph on $[11] \cup A$ given in Figure $2$, which includes precisely the edges $(i, a_{i})$ for $1 \leq i \leq 11$.  
\begin{center}
\includegraphics[scale=.5]{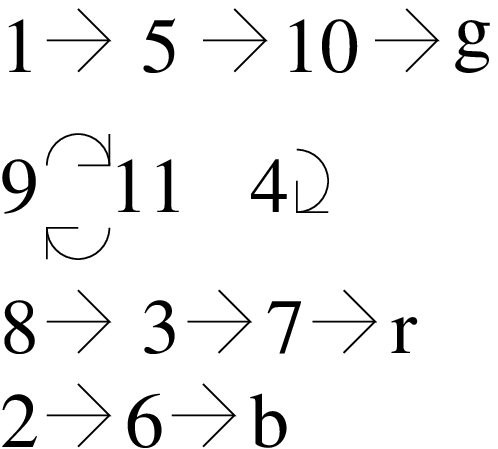}
\\Figure 2
\end{center}
Then $B^{g} = \left\{1, 4, 5, 9, 10, 11 \right\}$ (note that $1$ is connected to $g$, so $4, 9, 11$ are included in $B^{g}$), $B^{r} = \left\{3, 8, 7 \right\}$, and $B^{b} = \left\{2, 6 \right\}$, and $c_{1}^{r} = 3$, $c_{2}^{r} = 8$, $c_{3}^{r} = 7$; $c_{1}^{b} = 2$, $c_{2}^{b} = 6$.  The longest path containing $r$ is $8, 3, 7, r$ and the longest path containing $b$ is $2, 6, b$, so define a permutation $\pi \in S_{11}$ by $\pi = 5, 2, 8, 4, 10, 6, 3, 7, 11, 1, 9 = (1, 5, 10) (2) (3, 8, 7) (4) (6) (9, 11)$.  Define a coloring $\mu$ on the cycles of $\pi$ by coloring $(1, 5, 10)$, $(4)$, $(9, 11)$ using $g$; $(3, 8, 7)$ using $r$ and $(2)$, $(6)$ using $b$.  Then $(\pi, \mu)$ is a cycle-colored permutation of $[11]$.  
\section{The Generating Formula for Signed Stirling Numbers}
To prove \[
\sum_{\pi \in S_{n}} (-1)^{n-k(\pi)}x^{k(\pi)} = (x-n+1)(x-n+2) \cdots x
\] 
it is equivalent to prove
\[\sum_{\pi \in S_{n} \textrm{, } n-k(\pi) \textrm{ even}}x^{k(\pi)}= \sum_{\pi \in S_{n} \textrm{, } n-k(\pi) \textrm{ odd}}x^{k(\pi)}+ (x-n+1)(x-n+2) \cdots x \]
Because there are precisely $(x-n+1)(x-n+2) \cdots x$ cycle-colored permutations $(e_{n}, \mu)$ where $e_{n} \in S_{n}$ is the identity permutation on $[n]$ and $\mu$ uses precisely $n$ colors, it suffices to find a bijection \[\phi: \left\{(\pi, \mu): n-k(\pi) \textrm{even} \right\} \leftrightarrow \left\{(e_{n}, \mu): \mu \textrm{ uses } n \textrm { colors }\right\} \cup \left\{(\pi, \mu): n-k(\pi) \textrm{odd}\right\}\]
or, alternately, an involution $\phi$ on the set of cycle-colored permutations of $[n]$ such that $\phi((\pi, \mu)) = (\pi, \mu)$ if $\pi = e_{n}$ and $\mu$ uses precisely $n$ colors, and for all other $(\pi, \mu)$, if $\phi((\pi, \mu)) = (\tilde{\pi}, \tilde{\mu})$, $k(\pi) = k(\tilde{\pi}) \pm 1$.  

Define an automorphism $\phi$ on the set of cycle-colored permutations as follows:
\begin{enumerate}
\item	Given a cycle-colored permutation $(\pi, \mu)$, let $(i, j)$ with $1 \leq i < j \leq n$ be a \emph{relation} of $(\pi, \mu)$ if and only if the (not necessarily distinct) cycles containing $i$ and $j$ in $\pi$ have the same color under $\mu$, with $R = \left\{(i, j): (i, j) \textrm{ is a relation of } (\pi, \mu) \right\}$ the set of relations of $(\pi, \mu)$. 
\item  If $R = \emptyset$, let $\phi((\pi, \mu)) = (\pi, \mu)$   
\item If $R \neq \emptyset$, let $(i, j)$ be the minimal element of $R$ under the lexicographic ordering.  Define a permutation $\tilde{\pi}$ by $\tilde{\pi} = (ij) \circ \pi$, where $(ij)$ is the transposition sending $i$ to $j$ and vice-versa, and $\circ$ denotes multiplication in $S_{n}$.  
\item Define a coloring $\tilde{\mu}$ of the cycles of $\tilde{\pi}$ as follows: if $k \in [n]$ is in a cycle of $\pi$ colored using the color corresponding to $s$ under $\mu$, the cycle containing $k$ in $\tilde{\pi}$ is colored using the color corresponding to $s$ under $\tilde{\mu}$.  Let $\phi((\pi, \mu)) = (\tilde{\pi}, \tilde{\mu})$.     
\end{enumerate}
Note that, if $\pi = e_{n}$ and $\mu$ uses precisely $n$ colors, $R = \emptyset$, and $\phi((\pi, \mu)) = (\pi, \mu)$. For all other $(\pi, \mu)$, if $(i, j)$ is the minimal relation of $(\pi, \mu)$ and $\phi((\pi, \mu)) = (\tilde{\pi}, \tilde{\mu})$, either $i$ and $j$ are in the same cycle of $\pi$, in which case multiplying $\pi$ by $(ij)$ splits the cycle containing $i$ and $j$ into two cycles, and $k(\pi) = k(\tilde{\pi})-1$, or $i$ and $j$ are in distinct cycles of $\pi$, in which case multiplying $\pi$ by $(ij)$ joins the distinct cycle containing $i$ and the distinct cycle containing $j$ into one cycle, and $k(\pi) = k(\tilde{\pi})+1$.  Note also that, for all $k \in [n]$, the cycle containing $k$ in $\pi$ and the cycle containing $k$ in $\tilde{\pi}$ are colored the same way by, respectively, $\mu$ and $\tilde{\mu}$.  In particular, $(\pi, \mu)$ and $(\tilde{\pi}, \tilde{\mu})$ have the same set of relations and so the same minimal relation $(i, j)$, and $\phi$ is an involution.  This suffices to prove the above generating function of signed Stirling numbers.  
\subsection{Example} Let $n=5$ and $(\pi, \mu) = ((1) (2, 4, 3) (5), \mu)$ where $\mu$ colors $(1)$, $(5)$ using $r$ and $(2, 4, 3)$ using $b$.  Then the set of relations $R$ of $(\pi, \mu)$ is $\left\{ (1, 5), (2, 3), (2, 4), (3, 4) \right\}$, and the minimal element of $R$ under the lexicographic ordering is $(1, 5)$.  Define $\tilde{\pi} = (1, 5) (2, 4, 3)$, and $\tilde{\mu}$ as a coloring of the cycles of $\tilde{\pi}$, where $(1, 5)$ is colored using $r$ and $(2, 4, 3)$ is colored using $b$.  Then $(\tilde{\pi}, \tilde{\mu})$ is a cycle-colored permutation of $[5]$ with relations $\left\{ (1, 5), (2, 3), (2, 4), (3, 4) \right\}$ and $\phi((\pi, \mu)) = (\tilde{\pi}, \tilde{\mu})$, $\phi((\tilde{\pi}, \tilde{\mu})) = (\pi, \mu)$.  Note that $k(\pi) = 3$ and $k(\tilde{\pi}) = 2$.  
\section{Acknowledgements}
The author is indebted to Robin Pemantle for recognizing his color-based proof of the generating formula for signed Stirling numbers, to Herbert Wilf for suggesting he attempt to find a color-based proof of the generating formula for unsigned Stirling numbers; as well as to Janet Beissinger, who originated ~\cite{jbcolor} the view that the appearance of $x^{n}$ in a generating function could signal a coloring in $x$ colors.  
\bibliography{mybib}{}
\bibliographystyle{plain}
\end{document}